\documentclass[12pt]{amsart}
\usepackage{amscd}
\usepackage{verbatim}
\usepackage{amssymb, amsmath, amsthm, amscd,ifthen}
\usepackage[dvips]{graphics}
\usepackage[cp866]{inputenc}
\usepackage{graphicx}
\usepackage{epsfig}
\usepackage{color}

\usepackage{amsmath,amssymb,amscd,}
\usepackage[all]{xy}

\usepackage[dvips]{graphics}
\usepackage[cp866]{inputenc}
\usepackage{graphicx}
\usepackage{epsfig}
\usepackage[hang]{subfigure}

\theoremstyle{plain}
\newtheorem{theorem}{Theorem}
\newtheorem{lemma}{Lemma}

\theoremstyle{definition}

\theoremstyle{plain}
\newtoks\thehProclaim
\newtheorem*{Proclaim}{\the\thehProclaim}

\theoremstyle{definition}
\newtoks{\thehRemark}
\newtheorem*{Remark}{\the\thehRemark}

\renewcommand{\leq}{\leqslant}
\renewcommand{\geq}{\geqslant}

\begin{document}

\title{Extremal areas of polygons with fixed perimeter}

\author{Giorgi Khimshiashvili, Gaiane Panina, Dirk Siersma}

\address{G. Khimshiashvili: Ilia State University, Tbilisi; gogikhim@yahoo.com; G. Panina: St. Petersburg Department of Steklov Mathematical Institute; St. Petersburg State University;  gaiane-panina@rambler.ru;  D. Siersma: Utrecht University, Department of Mathematics; d.siersma@uu.nl}



\begin{abstract}
We consider the configuration space of planar $n$-gons with fixed perimeter, which is diffeomorphic to
the complex projective space $\mathbb{C}P^{n-2}$. The oriented area function
has the minimal number of critical points on
the configuration space. We describe its critical points (these are \textit{regular  stars}) and compute their indices when they are Morse.

\end{abstract}

\maketitle

\section{Introduction}

One of the results on isoperimetric problem discussed in the classical treatise of Legendre \cite{L}
states that the regular $n$-gon has the maximal area among all $n$-gons with
the fixed perimeter. For the historic development of this problem we refer to \cite{isoper}. The aim of present paper is to elaborate upon this classical result
by placing it in the context of Morse theory on a naturally associated configuration space.
To this end we follow the paradigms used in \cite{PanKhi}, \cite{panzh} and begin with several
definitions and recollections.

An $n$\textit{-gon} is an $n$-tuple of points  $(p_1,...,p_n)\in (\mathbb{R}^2)^n$, some of which may coincide.
Its \textit{perimeter} is (as usual)  $$\mathcal{P}(p_1,...,p_n)=|p_1p_2|+|p_2p_3|+...+|p_np_1|.$$

The configuration space $\mathcal{C}_n$ considered in the sequel is defined as the space of all polygons
(modulo rotations and translations) whose perimeter equals $1$  (one can fix any other positive number).

The \textit{oriented area} of a polygon with vertices $p_i=(x_i,y_i)$ is defined as
$$2A=x_1y_2-x_2y_1+...+x_ny_1-x_1y_n.$$

Oriented area as a Morse function has been studied in various settings:  for configuration spaces of
flexible polygons (those with side lengths fixed) in $\mathbb{R}^2$ and $\mathbb{R}^3$, see  \cite{PanGordTep},
\cite{PanKhi},  \cite{panzh},\cite{zhu}. We shall use some of the previous results in the new setting of
this paper.

\medskip

Let  $\sigma$ be a cyclic renumbering: given a polygon $P=(p_1,...,p_n)$, $$\sigma(p_1,...,p_n)=(p_2,p_3,...,p_n,p_1).$$
In other words, we have an action of $\mathbb{Z}_n$ on $\mathcal{C}_n$  which renumbers the vertices of a polygon.

\begin{figure}[h]
\centering \includegraphics[width=10 cm]{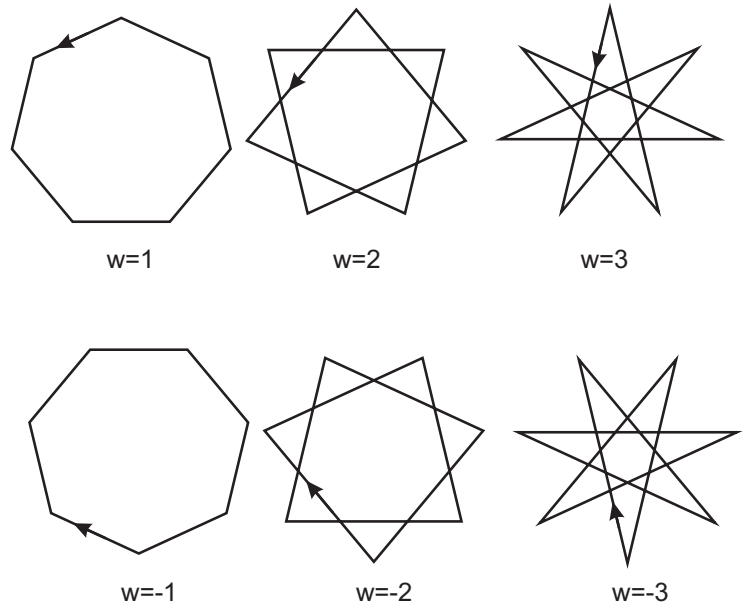}
\caption{Regular stars for $n=7$ with their winding numbers.}\label{Star7}
\end{figure}

\begin{figure}[h]
\centering \includegraphics[width=14 cm]{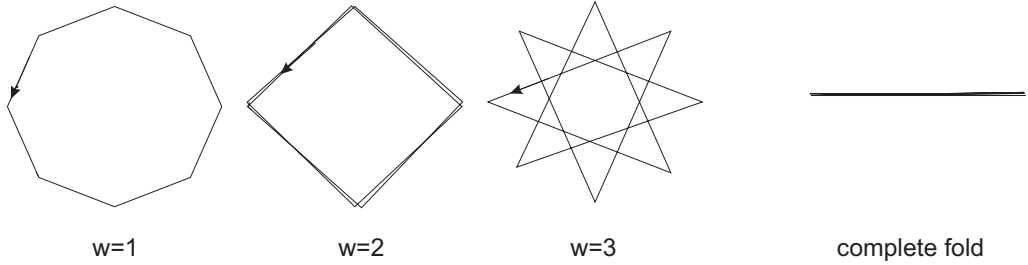}
\caption{Regular stars for $n=8$ with positive winding numbers.}\label{Star8}
\end{figure}

A \textit{regular star} is an equilateral $n$-gon  such that $\sigma (p_1,..., p_n)=(p_1,..., p_n)$, see Fig. \ref{Star7}, \ref{Star8}.

A \textit{complete fold} is a regular star with $p_i=p_{i+2}$. It exists for even $n$ only.

A regular star which is not a complete fold is uniquely defined by its winding number $w$ with respect to the center.

We can now formulate the main result of the paper.
\begin{theorem}  \begin{enumerate}
                   \item The  space $\mathcal{C}_n$  is homeomorphic to $\mathbb{C}P^{n-2}$. Therefore we consider it as a smooth manifold, keeping
                   in mind the smooth structure of the projective space,
                   \item Smooth critical points of  function $A$ on $\mathcal{C}_n$ are regular stars and complete folds only,
									
\item The non-smooth points of $A$ are Lipschitz-regular points of $A$,
\item The function $A$ has the minimal number of critical points on $\mathbb{C}P^{n-2}$,
                   \item The stars with the maximal winding numbers are non-degenerate critical points of function $A$,
                   \item Under assumption that all regular stars are non-degenerate critical points, the Morse index is:
                   $$M(P)=\left\{
                            \begin{array}{ll}
                              2w(P)-2, & \hbox{if $w(P)<0$;} \\
                             2n-2w(P)-2, & \hbox{if $w(P)>0$;} \\
                              n-2, & \hbox{if $P$ is a complete fold.}
                            \end{array}
                          \right.
                   $$
                 \end{enumerate}
\end{theorem}

As an illustration we give below drawings of some regular stars together with their winding numbers.

It is worthy mentioning that a scholarly example of an exact Morse function on complex projective space is $$F(u_1:...:u_{n-1})=\frac{a_1|u_1|^2+...+a_n|u_{n-1}|^2}{|u_1|^2+...+|u_{n-1}|^2}$$
for distinct real numbers $a_i$.  We will discuss more on this in the last section; let us now only mention that
an obvious difference with the area function  is that cyclic permutation
of indices preserves critical points of area and cyclically permutes critical points of $F$.

\medskip

\section{Proof  of Theorem 1}

There are three (well-known) statements to be used in the proof:

\medskip

\textbf{Statement A.} Assume a smooth manifold $C'$ is a codimension $k$ submanifold of a smooth  manifold $C$.
Let $q \in C'$ be a Morse point of some $f: C\rightarrow \mathbb{R}$. If $q$ is also a Morse point of the restriction
$f|_{C'}$, then for the Morse indices of $f$ related to $C$ and $C'$, we have:
$$M_{C'}(q)+k\geq M_C(q)\geq M_{C'}(q).$$

The proof easily follows from the Sylvester's formula for the index of symmetric matrix.
\medskip

\textbf{Statement B.} Assume a smooth manifold $C$ is a Cartesian product: $C=C_1\times C_2$, and the function
$f: C\rightarrow \mathbb{R}$ is the sum $f=f_1+f_2$, where $f_{1,2}: C_{1,2}\rightarrow \mathbb{R}$.
Then $q=(q_1,q_2)\in C$ is a Morse point if and only if $q_1$ and $q_2$ are Morse points of $f_1$ and $f_2$ respectively. In this case the Morse
indices sum up:
$$M(q)_C=M_{C_1}(q_1)+M_{C_2}(q_2).$$
This follows by noticing that, in this case, the Hessian matrix has a block form.

\bigskip

\textbf{Statement C.}  Assume that  $q$ is a point in  a smooth manifold $M^n$  which is an isolated critical point of a smooth function $f$.
Assume that there exist two submanifolds of complementary dimensions $N^k \subset M^n$ and $L^{n-k} \subset M^n$  intersecting transversally at $q$.
Assume also that $q$ is a non-degenerate local maximum  of  the restriction $f|_{N^k}$
and a non-degenerate local minimum  of  the restriction $f|_{L^{n-k}}$. Then $q$ is a Morse point of $f$ whose Morse index is $k$.

The proof again follows from the Sylvester's formula for the index of symmetric matrix.
\medskip

\textbf{The claim (1)}  of Theorem 1 was probably known before.  Assume we have a polygon $P\in \mathcal{C}_n$. Identifying its ambient space
$\mathbb{R}^2$  with $\mathbb{C}$, let us interpret its edges $p_1p_2$,...,$p_{n-2}p_{n-1}$  as complex numbers $u_1,...,u_{n-1}$.  This $(n-1)$-tuple
is never identical zero, so it gives a point $(u_1:...:u_{n-1}) \in \mathbb{C}P^{n-2}$. Conversely, each point  in the projective space yields a set of mutually homothetic polygons with a unique representative whose perimeter equals $1$.

\medskip

\textbf{The claim (2)} appeared already in \cite{Leger}  (yet unpublished). The proof follows from the two observations.

Firstly, it is easy to show that a critical point exhibits an equilateral polygon. Indeed, if this not the case, take two neighbour edges of different lengths, say, $l_1<l_2$, with $l_1=|p_1p_2|, \  l_2=|p_2p_3|$, freeze the rest of the polygon, and  move the point $p_2$  in such a way that $l_1+l_2$ remains fixed. The first differential of the area is non-zero since the point $l_2$ travels along an ellips, and the tangent line non-parallel to $p_1p_3$.

Secondly, if we freeze all the edge lengths of the polygon, we know from \cite{PanKhi} that being critical means being cyclic.
Therefore, a critical polygon is an equilateral polygon which is inscribed in a circle.

It remains to consider those having coinciding edges going in opposite direction, that is, with $p_i=p_{i+2}$.
These are complete folds.
The claim (2) is proven.

\bigskip

{\textbf{The claim (3)} is missing in \cite{Leger}.  The configuration space $\mathcal{C}_n$  inherits the smooth structure from $\mathbb{C}P^{n-2}$, but the function $A$ is not everywhere continuously differentiable, since the perimeter involves square roots. Non-smooth points are the configurations with two (or more) colliding consecutive vertices; they form an  arrangement of hyperplanes in $\mathbb{C}P^{n-2}$. To make sure that non-smooth points of $A$ "behave like regular ones" in Morse-theoretic sense   we will use the concept of Clarke subdifferential $\delta_c$ from non-smooth analysis. This subdifferential is related to the limits of all gradients in smooth points in a neighborhood of singular points. It was introduced by Clarke in \cite{clarke}. The generalized gradient  $\delta_c f$ of a locally Lipschiz function $f$ has values in the set of non-empty compact convex sets, satisfies the sum, product, quotient and chain rules with set addition and scalar multiplication. We use especially the reference \cite{APS}, which contains also the definition of critical point and {the} regular interval theorem.

In the homogeneous coordinates  $(u_1:...:u_{n-1}) \in \mathbb{C}P^{n-2}$ {the} signed area function  is  $\frac{A}{\mathcal{P}^2}$, where A is the area on $ \mathbb{C}^{n-1}$  and $\mathcal{P}$ the perimeter. More precisely:
$$ A_{\mathbb{C}P^{n-2}} (u_1:...:u_{n-1})  =
\frac{\sum_{1 \le i < j \le n-1} (\overline{u_i} u_j - \overline{u_j} u_i)}
{(|u_1|+ \cdots + |u_{n-1}| +  |u_1+ \cdots + u_{n-1}|)^2 } $$

Fix a  non-smooth point of $P$ and  let for this point  the first coordinate $u_1 = x_1 + i y_1=0$  and $w := u _2 + \cdots u_{n-1} \ne 0$.  Take now a local chart and
consider the real partial  subderivative with respect to $(x_1,y_1)$.
We use the following facts:
\begin{itemize}
\item[(a)]
$A'(0) = i\cdot w$  (usual derivative). This is $w$ rotated by $\frac{\pi}{2}$,
\item[(b)]
$ \delta_c  \sqrt{x_1^2 + y_1^2} = \mathbb D$, where  $\mathbb D$ is the unit disc in $\mathbb R^2$,
\item[(c)]
$\delta_c\mathcal{ P} (0) =\mathcal{P}'(0)  =  \frac{w}{|w|}  + \mathbb D$  (Clarke derivative via sum rule).
\end{itemize}
Next we us the quotient rule:
$$\Big(\frac{A}{\mathcal{P}^2}\Big)'=\frac{A'\cdot \mathcal{P}-2A \cdot \mathcal{P}'}{\mathcal{P}^3} =\mathcal{ P}_0 i \cdot w - 2 A_0  \frac{w}{|w|} - 2 A_0  \mathbb D.$$

{According} to \cite{APS} a point $q$ is called a {\it  critical point} for a locally Lipschitz function $f: M \to \mathbb R$ on a smooth manifold $M$  iff $O \in \delta_c(f) (q)$.   Proposition 1.2 of  \cite{APS} extends for Clarke regular points  the  first Morse lemma (regular interval theorem) to  locally Lipschitz functions.

We conclude {the} proof of (3) by showing that
a non-smooth point of $\mathcal{P}$  is never critical with respect to $\frac{A}{\mathcal{P}^2}$. It is sufficient to show this for the partial subderivative with respect to $(x_1,y_1)$.

The condition for critical point in this case  is :
$$  \mathcal{P}_0 i \cdot w - 2 A_0  \frac{w}{|w|}  \in  2 A_0  \mathbb D $$
By computing norms on both sides we get:
$\mathcal{P}_0^2 |w|^2 + 4 A_0^2  \leq 4 A_0^2$   and this is never the case. The claim (3) is proven.
\bigskip

{\textbf{The claim (4)}

We conclude that the number of critical points is $n-1$. Since the sum of the Betti-numbers and the LS-category of $\mathbb{C}P^{n-2}$ are also equal to $n-1$, the number of critical points is minimal.

\bigskip

{\textbf{The claims (5) and (6)}.}

Each of regular stars (which is not a complete fold) is uniquely defined by its winding number. The latter ranges from
$1$ to $[\frac{n-1}{2}]$  and from $-1$ to $-[\frac{n-1}{2}]$.

Therefore the Morse function is exact, and the Morse indices of the critical points are $0,2,4,...,2n-4$;
it remains to understand which star has which index.

\medskip

Denote by $S(n,w)$ the $n$-gonal regular star  with  index $w$.

First observe that each of the regular stars  (with the winding number $w$) has its symmetric image whose winding number is $-w$.

The only exception is the complete fold which is symmetric to itself.

\subsection*{Splitting construction. }

Take a neighborhood of a regular star  $P$ with a positive winding number $w$ (which is not a complete fold) and consider two transversally intersecting submanifolds of $\mathcal{C}_n$:
The first is the configuration space of equilateral polygons  $EQUILAT$.  That is, it consists of all polygons whose edge lengths are equal. Its dimension is $n-3$, and it is known from \cite{PanKhi}:

\begin{lemma}\label{FlexIndex} The Morse index of $P$ related to $EQUILAT$ is $n-1-2w$.\qed
\end{lemma}

Therefore, for $S(n,[\frac{n-1}{2}])$  the Morse index of $P$ related to this manifold is $0$. Moreover, from \cite{PanGordTep}  it is known that $P$ is a non-degenerate minimum.

\bigskip

The other manifold $CYCL$ is the subspace of all the polygons that are cyclic. Its dimension is $n-1$.

\begin{lemma}\label{SubmanifoldOfCyck}Let $P$ be a regular star which is not a complete fold with $w>0$.
Then $P$ is a non-degenerate local maximum on $CYCL$. So its Morse index of $P$ related to $CYCL$ is $n-1$.
\end{lemma}
Proof.  The problem is equivalent to the following setting:  Take the space of all $n$-gons inscribed in the unit circle.   Consider the function $A_{\mathbb{C}P^{n-2}}=\frac{A}{\mathcal{P}^2}$. It is well defined in a neighborhood of a regular star.
It suffices to prove that the regular star with $w>0$ is its local maximum.

Introduce local coordinates  by setting $x_i$  be the shift of the $i$-th vertex of the star, and   $t_i=x_{i+1}-x_{i}$. Note that $\sum_{i=1}^n x_i = 2 \pi $  and $\sum_{i=1}^n t_i= 0$.
We will compute  the $2$-jets $j_2$ (Taylor expansions up to the second degree). This is sufficient for finding the Morse indices.
Computations show that  in these  coordinates, $j_2\widetilde{A}$  is  a  negative definite
quadratic form.
In more detail,
$$j_2\mathcal{P}(t_1,...,t_n)= \sqrt{(2-2cos \ \alpha)} \cdot (n - 1/8 \sum_{i=1}^{n} t_i^2),$$
$$j_2A(t_1,...,t_n)=\frac{sin \ \alpha}{2} \cdot (2n - \sum_{i=1}^{n} t_i^2) ,$$
and
$$j_2(2A/\mathcal{P}^2)(t_1,...,t_n) = \frac{1+\cos \ \alpha}{8 n^2 \sin \ \alpha}\cdot (4n - \sum_{i=1}^{n} t_i^2), $$
where $\alpha=\frac{2w\pi}{n}$. Note that due to $w > 0$ the first factor is always positive and therefore the quadratic part is negative definite.

The  Lemma is proven.

Assume $n$ is odd.
Now application of Statement C proves Claims (5) and (6) of the Theorem 1 for $S(n,[\frac{n-1}{2}])$.
By symmetry, Claims (5) and (6) are also proven for $S(n,-[\frac{n-1}{2}])$.

\subsection*{ Proof of {Claim (6)} for the remaining cases}

We now assume that all critical points are non-degenerate.
 Two mutually symmetric stars have indices $i$ and $n-2-i$,
so it suffices to prove (6) for one of them. Besides, by symmetry reason we immediately conclude  that  the claim (6) holds for complete
folds

\medskip
First prove the statement for odd $n$. Take all the stars
$S(n,w)$  with $\frac{n-1}{2}>w>0$.  By Lemma \ref{SubmanifoldOfCyck}, $M(S(n,w))\geq  n-1$.
So their indices taken together are $2n-4,2n-6,...,n+1$.

Now let us take all the stars
$S(n,w)$  with $\frac{n-1}{2}<w<0$.  By symmetry, their indices taken together are $0,2,4,...,n-5$.
On the other hand, by Lemma \ref{FlexIndex}, their indices relative $EQUILAT$ give exactly the same set.
This relative index never exceeds the index rel $\mathcal{C}_n$, so Claim (6)  for $n$ odd is proven.

\medskip

Next prove the statement for even $n$. By a way of analogy consider all stars
$S(n,w)$  with $w>0$.  By Lemma \ref{SubmanifoldOfCyck}, $M(S(n,w))\geq  n-1$.
So their indices taken together are $2n-4,2n-6,...,n$.

Now let us take all the stars
$S(n,w)$  with $w<0$.  By symmetry, their indices taken together are $0,2,4,...,n-4$.
On the other hand, by Lemma \ref{FlexIndex}, their indices relative $EQUILAT$ give exactly the same set.
Since relative index never exceeds the index rel $\mathcal{C}_n$, Claim (6) is proven.
We have now completed the proof of Theorem 1. \qed

\section{Concluding remarks}

The study of Morse functions on  $\mathcal{C}_n=\mathbb{C}P^{n-2}$ can be related to its Veronese type embedding into the space HMat of Hermitian matrices. Its image is known (see \cite{Tight}) to be a taut embedding. It follows that  the restriction of almost every $\mathbb{R}$-linear function on HMat to the image of $\mathcal{C}_n=\mathbb{C}P^{n-2}$ is an exact Morse function, so  it has the minimal number of critical points. The function $F$ mentioned in the introduction is of this type. Our function $A$ is not of this type, since in the coordinates of  $\mathcal{C}_n=\mathbb{C}P^{n-2}$ is the
quotient of area and the square of the perimeter
$$ A_{\mathbb{C}P^{n-2}} (u_1:...:u_{n-1})  =
\frac{\sum_{1 \le i < j \le n-1} (\overline{u_i} u_j - \overline{u_j} u_i)}
{(|u_1|+ \cdots + |u_{n-1}| +  |u_1+ \cdots + u_{n-1}|)^2 } $$

The numerator extends to a linear function, but the denominator contains square roots of the coordinates. But it has still the good property that it has the minimal number of critical points. A more general question is how the geometry and  metric of projective space can be used for deeper understanding.

\medskip
It seems  also worthy to add that Theorem 1 suggests several further developments in the same spirit.

Firstly, one can obtain similar descriptions of extremals of area under different constraints. The most obvious option is to fix the sum of certain powers of side lengths. In other words, one looks for extremals of area under condition that $$\sum (a_i)^p = const,$$ where $p$ is a fixed positive number. For symmetry reasons, regular stars and complete folds are critical in this case too. However it is not obvious that there are no other extremals since the geometric reasoning used for fixed perimeter is not directly applicable. For $p=2,$ a rigorous analytic proof of this fact is given in \cite{Leger} without discussing non-degeneracy and Morse indices of stars.
Perhaps this case fits better in the framework of Veronese embedding.
 Moreover, it is not difficult to verify that the same method of proof is applicable for any positive $p$.

Finally, in line with the general idea of duality in the calculus of variations one may consider a dual problem by fixing the values of area and looking for extremals of perimeter. Some conclusions and lines of research are immediate. For example, it is easy to verify that generically the critical points of $A$  with $\mathcal{P}$ fixed are the same as the critical points of $\mathcal{P}$  with $A$ fixed. One can also derive the Morse index formula for the dual problem from our Theorem 1.

\section*{Acknowledgments}
It is our pleasure to acknowledge the  hospitality and excellent working conditions of
CIRM, Luminy, where this paper was initiated as a "Research in Pairs" project. This work is supported by the RFBR grant 17-01-00128.

\end{document}